\documentclass{article}

\usepackage{amsmath}
\usepackage{amsfonts}
\usepackage{amssymb}
\usepackage{eucal}

\newtheorem{thm}{Theorem}[section]
\newtheorem{cor}[thm]{Corollary}

\newtheorem{prop}[thm]{Proposition}
\newtheorem{rem}[thm]{Remark}
\newtheorem{ex}{Example}
\numberwithin{equation}{section}
\newenvironment{pr}{\noindent{\bf Proof }}{\space\hfill$\Box$\vskip2mm}

\begin{document}
\title{Slant curves in 3-dimensional normal almost paracontact metric manifolds}
\author{Joanna We{\l}yczko}






\date{\today}

\maketitle

\begin{abstract}
The presented paper is devoted to study the curvature and torsion of slant Frenet curves in 3-dimensional normal almost paracontact metric manifolds. Moreover, in this class of manifolds, properties of non-Frenet slant curves (with null tangents or null normals) are obtained. The achieved results are illustrated by examples.
\end{abstract}

\section{Preliminaries}

Let $M$ be a $(2n+1)$-dimensional differentiable manifold. Let $\varphi$ be  a $(1,1)$-tensor filed,  $\xi$ a vector filed and $\eta$ a 1-form  on $M$. Then $(\varphi,\xi,\eta)$ is called an almost paracontac structure on $M$ if
\begin{itemize}
\item [(i)] $\eta(\xi)=1, \quad \varphi^2 =\mathop{\rm Id}-\eta\otimes\xi$, 
\item [(ii)] the tensor filed $\varphi$ induces an almost paracomplex structure on the distribution $\mathcal D=\mathop{ker}\eta$, that is, the eigendistributions $\mathcal D^+$, $\mathcal D^-$ corresponding to the eigenvalues $1$, $-1$ of $\varphi$,
respectively, have equal dimension $n$.
\end{itemize}
$M$ is said to be almost paracontact manifold if it is endowed with an almost paracontact structure (cf.\ \cite{BCM,E,KW,Z}).
 
Let $M$ be an almost paracontact manifold. $M$ will be called an almost paracontact metric manifold if it is additionally endowed with a pseudo-Riemannian metric $g$ of signature $(n+1,n)$ and such that
$$ 
g(\varphi X,\varphi Y)=\null-g(X,Y)+\eta(X)\eta(Y).
$$
For such a manifold, we additionally have $\eta(X)=g(X,\xi)$, $\eta(\xi)=1$, $\varphi\xi=0$, $\eta\circ\varphi=0$. Moreover, we can define a skew-symmetric 2-form $\varPhi$ by $\varPhi(X, Y)=g(X,\varphi Y)$, which is called the fundamental form corresponding to the structure. Note that $\eta\wedge\Phi$ is up to a constant factor the Riemannian volume element of $M$.

On an almost paracontact manifold, one defines the $(2,1)$-tensor filed $N^{(1)}$ by 
\begin{equation*}
\label{nijen1}
N^{(1)}(X,Y)=[\varphi,\,\varphi](X,Y)-2\,d\eta(X,Y)\xi,
\end{equation*}
where $[\varphi, \varphi]$ is the Nijenhuis torsion of $\varphi$ given by 
$$
[\varphi,\,\varphi](X,Y)=\varphi^2[X,Y]+[\varphi X,\varphi Y]-\varphi[\varphi X,Y]-\varphi[X, \varphi Y].
$$ 
If $N^{(1)}$ vanishes identically, then the almost paracontact manifold (structure) is said to be normal (cf.\ \cite{KW} and recent papers \cite{BCM,Z}). The normality condition says that the almost paracomplex structure $J$ defined on $M\times\mathbb R$ by 
$$
  J\Big(X,\lambda\frac{d}{dt}\Big)=\Big(\varphi X+\lambda\xi,\eta(X)\frac{d}{dt}\Big)
$$
is integrable (paracomplex).

Since that time, we will consider manifolds of dimension $3$. The following theorem presents conditions equivalent to the normality of $3$ dimensional almost paracontact manifold, which we will use later in the work. 
\begin{thm}{\rm\cite{JW}}
\label{twab}
For a 3-dimensional almost paracontact metric manifold $M$, the following three conditions are mutually equivalent \\
\indent
(a) $M$ is normal, \\
\indent
(b) there exist functions $\alpha,\beta$ on $M$ such that
\begin{equation}
\label{prof2}
  (\nabla_X \varphi)Y=\beta(g(X,Y)\xi-\eta(Y)X) +\alpha(g(\varphi X,Y)\xi-\eta(Y)\varphi X),
\end{equation}
\indent
(c) there exist functions $\alpha,\beta$ on $M$ such that
\begin{equation}
\label{c}
  \nabla_X\xi=\alpha(X-\eta(X)\xi)+\beta\varphi X.
\end{equation}
\end{thm}
\begin{cor}{\rm\cite{JW}}
The functions $\alpha,\beta$ realizing (\ref{prof2}) as well as (\ref{c}) are given by 
\begin{equation}
\label{alphabeta}
  2\alpha =\mathop{\rm Trace}\,\{X\rightarrow\nabla_X\xi\},\qquad 
  2\beta =\mathop{\rm Trace}\,\{X\rightarrow~\varphi\nabla_X\xi\}.
\end{equation}
\end{cor}

A 3-dimensional normal almost paracontact metric manifold is said to be 
\begin{itemize}
\item paracosymplectic manifold if $\alpha=\beta=0$, \cite{Dacko}
\item quasi-para-Sasakian if and only if $\alpha=0$ and $\beta\neq 0$, \cite{E,JW}
\item $\beta$-para-Sasakian if and only if $\alpha=0$ and $\beta\neq 0$ and $\beta$ is constant, in  particular, para-Sasakian if $\beta=-1$ (\cite{JW,Z})
\item $\alpha$-para-Kenmotsu if $\alpha\neq 0$ and $\alpha$ is constant and $\beta=0$. 
\end{itemize}

At the end of this section, let $\gamma\colon I\to M$, $I$ being an interval, be a curve in $M$ such that  $g(\dot\gamma,\dot\gamma)=\pm1$ or $g(\dot\gamma,\dot\gamma)=0$ and let $\nabla_{\dot\gamma}$ denote the covariant differentiation along $\gamma$. 

We say that the curve $\gamma$ is a slant curve in $3$-dimensional normal almost paracontact metric manifold,  if 
\begin{equation*}
\label{cos}
g(\dot\gamma,\xi)=\eta(\dot\gamma)=c {\rm\; and\; {\it c}\; is\; constant.}
\end{equation*}
In the particular case, if $c=0$ the curve $\gamma$ is called Legendre curve \cite{JW}.

\begin{rem}
In normal almost contact metric manifold with Riemannian metric $g$, the value of $g(\dot\gamma,\xi)$ satisfies $-1\leqslant g(\dot\gamma,\xi)\leqslant 1$, so that we can define the structural angel of $\gamma$ i.e the function $\mathop{\theta}:I\rightarrow[0,2\pi)$ given by 
$$
\cos \theta(t)=g(\dot\gamma(t),\xi)=\eta(\dot\gamma(t)).
$$
Then, the curve $\gamma$ is said to be a slant curve, ( or $\theta$-slant curve), if $\theta$ is a constant function, (see {\rm  \cite{cc, cc2, ccm}}). 
\end{rem}

Let $\gamma:I\rightarrow M$ be a slant curve in $M$ such that $g(\dot\gamma,\dot\gamma)=\varepsilon_1=\pm1$.
Let us consider tree vector fields $\dot\gamma, \varphi\dot\gamma, \xi$ for which, we have
$$
g(\dot\gamma,\dot\gamma)=\varepsilon_1,\; g(\dot\gamma,\xi)=c,\; g(\varphi\dot\gamma,\varphi\dot\gamma)=-\varepsilon_1+c^2,\; g(\xi,\xi)=1,\; g(\dot\gamma,\varphi\dot\gamma)=g(\xi,\varphi\dot\gamma)=0.
$$

This vector field are linearly independent and form the base of $TM$ if and only if $\varepsilon_1-c^2\neq 0$. In this case, we define  an orthonormal frame along $\gamma$ as follows:
\begin{equation}
\label{baza F}
F_1=\dot\gamma, \quad F_2=\frac{\varphi\dot\gamma}{\sqrt{|\varepsilon_1-c^2|}},\quad F_3=\frac{\xi-\varepsilon_1c\dot\gamma}{\sqrt{|\varepsilon_1-c^2|}},
\end{equation}
where $ g(F_1,F_1)=\varepsilon_1$, $g(F_2, F_2)=\mathop{sgn}{(-(\varepsilon_1-c^2))}=\upsilon$, $g(F_3, F_3)=-\varepsilon_1\upsilon$.

The frame $(F_1,F_2,F_3)$ will play an auxiliary role in the proofs of  Theorems \ref{kt} and \ref{nn}.
We additionally have
\begin{prop}
Let  $M$ be a 3-dimensional normal almost paracontact metric manifold.
If $\gamma\colon I\to M$ is a slant Frenet curve of osculating order 3 which is not an integral curve of $\xi$, then
\begin{eqnarray}
\label{f1}
\nabla_{\dot\gamma}F_1&=&\null\upsilon\delta\sqrt{|\varepsilon_1-c^2|}F_2-\varepsilon_1\alpha \sqrt{|\varepsilon_1-c^2|}F_3,\\
\label{f2}
\nabla_{\dot\gamma}F_2&=&\null-\varepsilon_1\delta\sqrt{|\varepsilon_1-c^2|}F_1+(\varepsilon_1\beta-\upsilon c\delta)F_3,\\
\label{f3}
\nabla_{\dot\gamma}F_3&=&-\varepsilon_1\upsilon\alpha\sqrt{|\varepsilon_1-c^2|}F_1+(\beta-\varepsilon_1\upsilon c \delta )F_2,
\end{eqnarray}
where  $\delta$ is a function defined by $\delta=\dfrac{g(\nabla_{\dot\gamma}{\dot\gamma},\varphi\dot\gamma)}{|\varepsilon_1-c^2|}$, and for simplicity we write $\alpha$, $\beta$ instead of the composed functions $\alpha\circ\gamma$, $\beta\circ\gamma$ with $\alpha$, $\beta$ being the same as in (\ref{alphabeta}).
\end{prop}

Separately, we discuss the case when $\dot\gamma, \varphi\dot\gamma, \xi$ are linearly dependent.
Let us note that  $\dot\gamma, \varphi\dot\gamma, \xi$ are linearly dependent if and only if
\begin{equation}
\label{zal}
\dot\gamma=c\xi \quad {\rm or} \quad \dot\gamma=c\xi\pm\varphi\dot\gamma.
\end{equation}

Indeed, if  $\dot\gamma, \varphi\dot\gamma, \xi$ are linearly dependent then  $\varepsilon_1=1$ and $c^2=1$. Hence $g(\varphi\dot\gamma,\varphi\dot\gamma)=0$ and either $\varphi\dot\gamma=0$ or $\varphi\dot\gamma$ is an isotropic vector field. 

If $\varphi\dot\gamma=0$ then $0=\varphi^2\dot\gamma=\dot\gamma-c\xi$ and $\dot\gamma=c\xi$. 

If $\varphi\dot\gamma$ is an isotropic vector field then,  we can write $\dot\gamma=a\xi+b\varphi\dot\gamma$, for some functions $a$ and $b$. Using $g(\dot\gamma,\dot\gamma)=1$ and $g(\dot\gamma,\xi)=c$, we find $\dot\gamma=c\xi+b\varphi\dot\gamma$. Now 
$$
\varphi\dot\gamma=b\varphi^2\dot\gamma=b(\dot\gamma-c\xi)=b(c\xi+b\varphi\dot\gamma-c\xi)=b^2\varphi\dot\gamma.
$$
Hence $b^2=1$ and $\dot\gamma=c\xi\pm\varphi\dot\gamma$.

\section{Frenet slant curves}

Let $M$ be a $3$-dimensional pseudo-Riemannian manifold with metric $g$. 

Let $\gamma\colon I\to M$, $I$ being an interval, be a curve in $M$. 

We say that $\gamma$ is a Frenet curve if $g(\dot\gamma,\dot\gamma)=\varepsilon_1$, $\varepsilon_1=\pm 1$ and one of the following three cases hold \\
\indent(a) $\gamma$ is of osculating order 1, i.e., $\nabla_{\dot\gamma}\dot\gamma=0$ (geodesics); \\
\indent(b) $\gamma$ is of osculating order 2, i.e., there exist two orthonormal vector fields $E_1(=\dot\gamma)$, $E_2$, $(g(E_2,E_2)=\varepsilon_2=\pm 1)$ and a positive function $\kappa$ (the curvature) along $\gamma$ such that 
$$
  \nabla_{\dot\gamma}E_1 = \kappa\varepsilon_2 E_2,\quad 
  \nabla_{\dot\gamma}E_2 = -\kappa\varepsilon_1 E_1;
$$
\indent(c) $\gamma$ is of osculating order 3, i.e., there exist three orthonormal vector fields $E_1(=\dot\gamma)$, $E_2$, $E_3$, ($g(E_2,E_2)=\varepsilon_2=\pm 1$, $g(E_3,E_3)=\varepsilon_3=\pm 1$) and two positive functions $\kappa$ (the curvature) and $\tau$ (the torsion) along $\gamma$ such that
$$
  \nabla_{\dot\gamma}E_1 = \kappa \varepsilon_2 E_2,\quad 
  \nabla_{\dot\gamma}E_2 = \null-\kappa\varepsilon_1 E_1+\tau \varepsilon_3E_3,\quad 
  \nabla_{\dot\gamma}E_3 = \null-\tau\varepsilon_2 E_2.
$$

\begin{prop}
\label{geodezyjne}
Let  $M$ be a 3-dimensional normal almost paracontact metric manifold.
If $\gamma\colon I\to M$ is a slant Frenet curve  in $M$, such that $g(\dot\gamma,\dot\gamma)=1$ and $g(\dot\gamma,\xi)=c=\pm1$ then $\gamma$ is a geodesics.
\end{prop}

\begin{pr}
In view of (\ref{zal}), we consider two opportunities. 

If  $\dot\gamma=c\xi$ then from (\ref{c}) $\nabla_{\dot\gamma}{\dot\gamma}=\nabla_{\xi}{\xi}=0$ and $\gamma$ is a geodesics.

If $\dot\gamma=c\xi\pm\varphi\dot\gamma$ using (\ref{c}) and (\ref{prof2}), we get 
\begin{eqnarray*}
\nabla_{\dot\gamma}{\dot\gamma}&=&c\nabla_{\dot\gamma}\xi\pm\nabla_{\dot\gamma}\varphi\dot\gamma=c\nabla_{\dot\gamma}\xi\pm(\nabla_{\dot\gamma}\varphi)\dot\gamma\pm\varphi\nabla_{\dot\gamma}{\dot\gamma}\\
&=&(\alpha\mp\beta)(c\dot\gamma-\xi\mp c\varphi\dot\gamma)\pm\varphi\nabla_{\dot\gamma}{\dot\gamma}=\null\pm\varphi\nabla_{\dot\gamma}{\dot\gamma}.
\end{eqnarray*}
Then $g(\nabla_{\dot\gamma}{\dot\gamma},\nabla_{\dot\gamma}{\dot\gamma})=\null\pm g(\nabla_{\dot\gamma}{\dot\gamma},\varphi\nabla_{\dot\gamma}{\dot\gamma})=0$ and $\gamma$ is a Frenet curve so that $\gamma$ is a geodesics.
\end{pr}

 Theorem \ref{kt} generalizes results obtained for the  Legendre curve  in \cite{JW}.

\begin{thm}
\label{kt}
Let  $M$ be a 3-dimensional normal almost paracontact metric manifold.
If $\gamma\colon I\to M$ is a slant Frenet curve of osculating order 3 in $M$, then its curvature and torsion are given by
\begin{eqnarray}
\label{kappa}
\kappa&=&\sqrt{|\varepsilon_1-c^2||\alpha^2-\varepsilon_1\delta^2|},\\
\label{tau}
\tau&=&\Big|\mathop{sgn}(1-\varepsilon_1 c^2)\beta+c\delta+\frac{\alpha\dot\delta-\dot\alpha\delta}{\alpha^2-\varepsilon_1\delta^2}\Big|.
\end{eqnarray}
where  $\delta$ is a function defined by $\delta=\dfrac{g(\nabla_{\dot\gamma}{\dot\gamma},\varphi\dot\gamma)}{|\varepsilon_1-c^2|}$, and for simplicity we write $\alpha$, $\beta$ instead of the composed functions $\alpha\circ\gamma$, $\beta\circ\gamma$ with $\alpha$, $\beta$ being the same as in (\ref{alphabeta}).
\end{thm}

\begin{pr}
Let $\gamma$ be a slant Frenet curve of osculating order 3 in $M$.
From Proposition \ref{geodezyjne}, we know, that for such a curve $\varepsilon_1-c^2\neq0$ and so that we can use the  frame $(F_1, F_2, F_3)$ from (\ref{baza F}) to calculate the curvature function and the torsion function of a curve  $\gamma$. 

We obtain the curvature function using (\ref{f1}) and calculating the length of $\nabla_{\dot\gamma}{\dot\gamma}$. 
We get 
$$\kappa^2\varepsilon_2=-\varepsilon_1\upsilon|\varepsilon_1-c^2|(\alpha^2-\varepsilon_1\delta^2)$$ 
so that the curvature is given by (\ref{kappa}) and  $\varepsilon_2=\mathop{sgn}(-\varepsilon_1\upsilon(\alpha^2-\varepsilon_1\delta^2))=\pm1$. 

Using (\ref{f1}), (\ref{f2}) and (\ref{f3}), we find
\begin{eqnarray*}
E_2=\frac{1}{\varepsilon_2\kappa}\nabla_{\dot\gamma}{E_1}=\frac{\varepsilon_2(\upsilon\delta F_2-\varepsilon_1\alpha F_3)}{\sqrt{|\alpha^2-\varepsilon_1\delta^2|}},
\end{eqnarray*}

Let us denote $p=\sqrt{|\alpha^2-\varepsilon_1\delta^2|}$ and calculate
\begin{eqnarray}
\label{fr2}
\nabla_{\dot\gamma}E_2&=&\varepsilon_2\upsilon\dot\gamma\Big(\frac{\delta}{p}\Big)F_2+\varepsilon_2\upsilon\Big(\frac{\delta}{p}\Big)\nabla_{\dot\gamma}F_2\\\nonumber
&&\null-\varepsilon_1\varepsilon_2\dot\gamma\Big(\frac{\alpha}{p}\Big)F_3-\varepsilon_1\varepsilon_2\Big(\frac{\alpha}{p}\Big)\nabla_{\dot\gamma}F_3\\\nonumber
&=&-\varepsilon_1\kappa E_1+\Big(\varepsilon_2\upsilon\dot\gamma\Big(\frac{\delta}{p}\Big)-\varepsilon_2\frac{\alpha\upsilon}{p}\Big(\varepsilon_1\upsilon\beta- c\delta\Big)\Big)F_2\\\nonumber
&&\null+\Big(-\varepsilon_1\varepsilon_2\dot\gamma\Big(\frac{\alpha}{p}\Big)+\varepsilon_2\frac{\delta}{p}\Big(\varepsilon_1\upsilon\beta-c\delta\Big)\Big)F_3\\\nonumber
&=&-\varepsilon_1\kappa E_1
+\varepsilon_2\Big(-\varepsilon_1\upsilon\beta+c\delta+\frac{\alpha\dot\delta-\dot\alpha\delta}{\alpha^2-\varepsilon_1\delta^2}\Big)\frac{\upsilon\alpha F_2-\delta F_3}{\sqrt{|\alpha^2-\varepsilon_1\delta^2|}}.
\end{eqnarray}
Hence from (\ref{fr2}), in view of Frenet equations, we have
\begin{eqnarray*}
\varepsilon_3\tau^2&=&\upsilon\mathop{sgn}(\alpha^2-\varepsilon_1\delta^2)\Big(-\varepsilon_1\upsilon\beta+c\delta+\frac{\alpha\dot\delta-\dot\alpha\delta}{\alpha^2-\varepsilon_1\delta^2}\Big)^2\\\nonumber
&=&-\varepsilon_1\varepsilon_2\Big(\mathop{sgn}(1-\varepsilon_1 c^2)\beta+c\delta+\frac{\alpha\dot\delta-\dot\alpha\delta}{\alpha^2-\varepsilon_1\delta^2}\Big)^2.
\end{eqnarray*}
This gives (\ref{tau}) and $\varepsilon_3=-\varepsilon_1\varepsilon_2$.
\end{pr}

Below, we present the functions of curvature and torsion of a slant Frenet curves of osculating order $3$   in some subclasses of 3-dimensional normal almost paracontact metric manifold.

\begin{cor}
Let  $M$ be a 3-dimensional manifold,
$\gamma\colon I\to M$ be a slant Frenet curve of osculating order 3 in $M$. 
 
 If $M$ is a paracosymplectic manifold then the curvature and torsion of $\gamma$ are given by
\begin{eqnarray*}
\kappa&=&\sqrt{|\varepsilon_1-c^2|}|\delta|,\quad
\tau=|c\delta|.
\end{eqnarray*}

If $M$ is a quasi-para-Sasakian manifold then the curvature and torsion of $\gamma$ are given by
\begin{eqnarray*}
\kappa&=&\sqrt{|\varepsilon_1-c^2|}|\delta|,\quad
\tau=\Big|\mathop{sgn}(1-\varepsilon_1 c^2)\beta+c\delta\Big|.
\end{eqnarray*}

If $M$ is a para-Sasakian manifold then the curvature and torsion of $\gamma$ are given by
\begin{eqnarray*}
\kappa&=&\sqrt{|\varepsilon_1-c^2|}|\delta|,\quad
\tau=\Big|\mathop{sgn}(1-\varepsilon_1 c^2)+c\delta\Big|.
\end{eqnarray*}

If $M$ is an $\alpha$-Kenmotsu manifold then the curvature and torsion of $\gamma$ are given by
\begin{eqnarray*}
\kappa&=&\sqrt{|\varepsilon_1-c^2||\alpha^2-\varepsilon_1\delta^2|},\quad
\tau=\Big|c\delta+\frac{\alpha\dot\delta}{\alpha^2-\varepsilon_1\delta^2}\Big|.
\end{eqnarray*}
Moreover $\delta$ is a function defined by $\delta=\dfrac{g(\nabla_{\dot\gamma}{\dot\gamma},\varphi\dot\gamma)}{|\varepsilon_1-c^2|}$, and for simplicity we write $\alpha$, $\beta$ instead of the composed functions $\alpha\circ\gamma$, $\beta\circ\gamma$ with $\alpha$, $\beta$ being the same as in (\ref{alphabeta}).
\end{cor}
\section{ Null slant Curves }

Next, we consider the case when $\gamma$ is a null curve, i.e has a null tangent ( $g(\dot\gamma,\dot\gamma)=0$).

We consider null curve  $\gamma$ which is not a geodesic, then $g(\nabla_{\dot\gamma}\dot\gamma, \nabla_{\dot\gamma}\dot\gamma)\neq0$. We take a parametrization of $\gamma $ such that $g(\nabla_{\dot\gamma}\dot\gamma, 
\nabla_{\dot\gamma}\dot\gamma)=1$ and define (see \cite{D,DB,DJ,KH} )
\begin{equation}
\label{def}
T=\dot\gamma,\quad N=\nabla_{\dot\gamma}T,\quad \tau=\frac{1}{2} g(\nabla_{\dot\gamma}N,\nabla_{\dot\gamma}N),\quad
W=-\nabla_{\dot\gamma}N-\tau T,
\end{equation}
where 
$$g(T,W)=g(N,N)=1 \quad {\rm and}\quad g(T,T)=g(T,N)=g(W,W)=g(W,N)=0.$$
Then, we obtain the Cartan equations for null curve $\gamma$.
\begin{eqnarray*}
\nabla_{\dot\gamma}T=N,\quad
\nabla_{\dot\gamma}W=\tau N.\quad
\nabla_{\dot\gamma}N=-\tau T- W.
\end{eqnarray*}

In \cite{JW}, we proved that every null Legendre curve in 3-dimensional normal almost paracontact metric manifold is a geodesic.

In this section, we consider a slant non-geodesic null curve $\gamma:I\rightarrow M$.

\begin{thm}
If $\gamma:I\rightarrow M$ is a slant non-geodesic null curve in a $3$-di\-men\-sio\-nal normal almost paracontact metric manifold then
\begin{eqnarray}
\label{tau null}
\tau&=&-\frac{\alpha^2c^2}{2}-\dot\alpha c\pm\beta-\frac{1}{2c^2},\\
\label{N}
N&=&\alpha c\dot\gamma\pm\frac{1}{c}\varphi\dot\gamma,\\
\label{W}
W&=&\frac{-\alpha^2c^2-1}{2c}\dot\gamma\mp\alpha\varphi\dot\gamma+\frac{1}{c^2}\xi,
\end{eqnarray}
where for simplicity we write $\alpha$, $\beta$ instead of the composed functions $\alpha\circ\gamma$, $\beta\circ\gamma$ with $\alpha$, $\beta$ being the same as in (\ref{alphabeta}).
\end{thm}

\begin{pr}
First, we find a vector field $N$. Let $\varphi\dot\gamma=p\dot\gamma+q N+r W$ for some functions $p,\,q,\,r$.
We find 
$$r=g(\varphi\dot\gamma,\dot\gamma)=0,\quad 
 g(\varphi\dot\gamma,\varphi\dot\gamma)=c^2= q^2,\quad g(\varphi\dot\gamma,\xi)=pc+q\alpha c^2=c(p\null\pm\alpha c^2)=0.
$$
Hence  $r=0,\; q=\null\pm c,\; p=\null\mp\alpha c^2$ and  $\varphi\dot\gamma=\null\mp\alpha c^2\dot\gamma\pm c N$, which leads to (\ref{N}).

Now, applying (\ref{prof2}) and (\ref{N}), we calculate
\begin{eqnarray}
\label{nabla N}
\nabla_{\dot\gamma}N&=&\dot\alpha c\dot\gamma+\alpha c\nabla_{\dot\gamma}{\dot\gamma}\pm\frac{1}{c}((\nabla_{\dot\gamma}\varphi)\dot\gamma+\varphi\nabla_{\dot\gamma}{\dot\gamma})\\\nonumber
&=&\Big(\dot\alpha ca+\alpha^2 c^2\mp\beta+\frac{1}{c}\Big)\dot\gamma\pm\alpha\varphi\dot\gamma-\frac{1}{c}\xi.
\end{eqnarray}
Next, in view of definition $\tau$ and $W$ in (\ref{def}), we use (\ref{nabla N}) and get (\ref{tau null}) and (\ref{W}).
\end{pr}

In the following corollary, we present the properties of slant null curves in some subclasses of 3-dimensional normal almost paracontact metric manifold.

\begin{cor}
Let  $M$ be a 3-dimensional  manifold,
 $\gamma\colon I\to M$ be a slant non-geodesic null curve in $M$. 
 
If $M$ is a paracosymplectic manifold then for  $\gamma$, we have
\begin{eqnarray*}
&&\tau=-\frac{1}{2c^2},\quad 
N=\null\pm\frac{1}{c}\varphi\dot\gamma,\quad 
W=\frac{1}{c^2}\xi.
\end{eqnarray*}

If $M$ is a quasi-para-Sasakian manifold then for  $\gamma$, we have
\begin{eqnarray*}
&&\tau=\null\pm\beta-\frac{1}{2c^2},\quad
N=\null\pm\frac{1}{c}\varphi\dot\gamma,\quad
W=\frac{1}{c^2}\xi.
\end{eqnarray*}

If $M$ is a para-Sasakian manifold then for  $\gamma$, we have
\begin{eqnarray*}
&&\tau=\null\pm1-\frac{1}{2c^2},\quad
N=\null\pm\frac{1}{c}\varphi\dot\gamma,\quad
W=\frac{1}{c^2}\xi.
\end{eqnarray*}

If $M$ is an $\alpha$-Kenmotsu manifold then for  $\gamma$, we have
\begin{eqnarray*}
&&\tau=-\frac{\alpha^2c^2}{2}-\frac{1}{2c^2},\quad
N=\alpha c\dot\gamma\pm\frac{1}{c}\varphi\dot\gamma,\quad
W=\frac{-\alpha^2c^2-1}{2c}\dot\gamma\mp\alpha\varphi\dot\gamma+\frac{1}{c^2}\xi.
\end{eqnarray*}.
\end{cor}

\section{Slant curves with null normal}

At last,  we investigate curves with null normals (for such curves in Minkowski spaces $\mathbb E^3_1$, see among others \cite{bon,IN,int}. We say $\gamma\colon I\to M$ is a curve with null normal if 
\begin{equation*}
\label{nulln}
  g(\dot\gamma,\dot\gamma)=\varepsilon_1=\pm 1,\quad 
  \nabla_{\dot\gamma}\dot\gamma \neq0, \quad 
  g(\nabla_{\dot\gamma}\dot\gamma,\nabla_{\dot\gamma}\dot\gamma)=0.
\end{equation*}

For such a curve, we define  
$$
T=\dot\gamma,\quad N=\nabla_{\dot\gamma}T,\quad W,
$$
where 
\begin{equation}
\label{cartan}
\begin{array}{l}
g(T,T)=\varepsilon_1, \quad g(N,W)=1,\\
 g(N,N)=g(W,W)=g(T,N)=g(T,W)=0.
\end{array}
\end{equation}
The signature of metric $g$ is $(-,+,+)$ so that $\varepsilon_1=1$.

We obtain the following Cartan equations for null normal curve $\gamma$. (In Minkowski 3-space see \cite{int})
\begin{eqnarray*}
\nabla_{\dot\gamma}T=N,\quad
\nabla_{\dot\gamma}N=\kappa N,\quad
\nabla_{\dot\gamma}W=- T-\kappa W.
\end{eqnarray*}

\begin{thm}
\label{nn}
If $\gamma:I\rightarrow M$ is a slant non-geodesic curve with null normal in a $3$-dimensional normal almost paracontact metric manifold, for which $c^2\neq1$ then $\alpha=\null\pm\dfrac{g(N,\varphi\dot\gamma)}{1-c^2}\neq0$,  and

\begin{eqnarray}
\label{N nor}
N&=&\null-\alpha(\xi-c\dot\gamma\pm\varphi\dot\gamma),\\
\label{kappa nor}
\kappa&=&\frac{\dot\alpha}{\alpha}\pm\beta+\alpha c,\\
\label{W nor}
W&=&\null\frac{-1}{2\alpha(1-c^2)}(\xi-c\dot\gamma\mp\varphi\dot\gamma),
\end{eqnarray}
where for simplicity we write $\alpha$, $\beta$ instead of the composed functions $\alpha\circ\gamma$, $\beta\circ\gamma$ with $\alpha$, $\beta$ being the same as in (\ref{alphabeta}).
\end{thm}


\begin{pr}
For a slant curve with null normal for which $c^2\neq1$, we can use an  orthonormal frame along $\gamma$ as in (\ref{baza F}).
We put 
\begin{equation}
\label{n pom}
N=\nabla_{\dot\gamma}{\dot\gamma}=a F_2+b F_3,
\end{equation}
for some function $a$ and $b$ and calculate $g(N,N)=a^2-b^2=0$. Hence  $a=\pm b$. 
Moreover, using (\ref{c}), we get
\begin{eqnarray}
\label{b}
b=g(N,F_3)&=&\null-\frac{1}{\sqrt{|1-c^2|}}g(\nabla_{\dot\gamma}\xi,\dot\gamma)\\\nonumber
&=&\null-\alpha\sqrt{|1-c^2|}.
\end{eqnarray} 
Hence $\alpha\neq0$ and putting (\ref{b}) to (\ref{n pom}), we obtain (\ref{N nor}).

To find the curvature, we use (\ref{c}), (\ref{N nor}) , (\ref{prof2}) and calculate
\begin{eqnarray*}
\nabla_{\dot\gamma} N&=&\null-\dot\alpha(\xi-c\dot\gamma\pm\varphi\dot\gamma)\\\nonumber
&&\null-\alpha(\nabla_{\dot\gamma}\xi-c\nabla_{\dot\gamma}{\dot\gamma}\pm(\nabla_{\dot\gamma}\varphi)\dot\gamma\pm\varphi\nabla_{\dot\gamma}{\dot\gamma})\\\nonumber
&=&(-\dot\alpha\mp\alpha\beta-\alpha^2 c)(\xi-c\dot\gamma\pm\varphi\dot\gamma)\\\nonumber
&=&\Big(\frac{\dot\alpha}{\alpha}\pm\beta+\alpha c\Big)N.
\end{eqnarray*}
In view of Cartan frame of $\gamma$, we obtain (\ref{kappa nor}).

At least, we put $W=d F_1+e F_2+f F_3$ for some functions $d,e,f$ and using (\ref{cartan}), we get (\ref{W nor}).
\end{pr}

As an immediate consequence of the above theorem, we obtain

\begin{thm}
Let $\gamma:I\rightarrow M$ be a Legendre curve with null normal in a 3-di\-men\-sio\-nal normal almost paracontact metric manifold $M$. Then $\alpha=\pm g(N,\varphi\dot\gamma)$, $\alpha\neq0$  and  we respectively have 
$$
T=\dot\gamma,\quad N=-\alpha(\xi\pm\varphi\dot\gamma), \quad W=\frac{-1}{2\alpha}(\xi\mp\varphi\dot\gamma)\quad {\it and}\quad \kappa=\left(\frac{\dot\alpha}{\alpha}\pm\beta\right), 
$$
where for simplicity we write $\alpha$, $\beta$ instead of the composed functions $\alpha\circ\gamma$, $\beta\circ\gamma$ with $\alpha$, $\beta$ being the same as in (\ref{alphabeta}).
\end{thm}

\begin{cor}
Let $\gamma:I\rightarrow M$ be a slant curve with null normal in a 3-dimensional $\alpha$-Kenmotsu manifold $M$, for which $c^2\neq1$. Then 
$\alpha=\null\pm g(N,\varphi\dot\gamma)\in\mathbb{R}$ and
\begin{eqnarray*}
&&\kappa=\alpha c,\quad
N=-\alpha(\xi-c\dot\gamma\pm\varphi\dot\gamma),\quad
W=\frac{-1}{2\alpha(1-c^2)}(\xi-c\dot\gamma\mp\varphi\dot\gamma).
\end{eqnarray*}
In the particular case when $\gamma$ is a Legendre curve, we have 
\begin{eqnarray*}
&&\kappa=0,\quad
N=-\alpha(\xi\pm\varphi\dot\gamma),\quad
W=\frac{-1}{2\alpha}(\xi\mp\varphi\dot\gamma).
\end{eqnarray*}
where for simplicity we write $\alpha$, $\beta$ instead of the composed functions $\alpha\circ\gamma$, $\beta\circ\gamma$ with $\alpha$, $\beta$ being the same as in (\ref{alphabeta}).
\end{cor}

\section{Examples illustrating theorems}

Legendre curves with  null normal that will appear in examples are  from my previous paper \cite{JW}.

\begin{ex}
{\rm
Let $\mathbb R^3$ be the Cartesian space and $(x,y,z)$ be the Cartesian coordinates in it. Define the standard almost paracontact structure $(\varphi,\xi,\eta)$ on $\mathbb R^3$ by 
\begin{equation}
\label{fiksieta}
  \varphi\partial_1=\partial_2-2x\partial_3,\quad
  \varphi\partial_2=\partial_1,\quad
  \varphi\partial_3=0,\quad
  \xi=\partial_3,\quad
  \eta=2xdy+dz,
\end{equation}
where $\partial_1 = \dfrac{\partial}{\partial x}$, $\partial_2 = \dfrac{\partial}{\partial y}$ and $\partial_3 = \dfrac{\partial}{\partial z}$. By certain direct calculations, one verifies that 
$$
  [\varphi,\varphi](\partial_i,\partial_j)-2d\eta(\partial_i,\partial_j)\xi=0,\quad 1\leqslant i<j\leqslant3,
$$
so that (\ref{nijen1}) is satisfied and the structure is normal. 

Suppose that $M=\mathbb R^2\times\mathbb R_+\subset\mathbb R^3$ and consider a normal almost paracontact metric structure on $M$ defined in the following way: $(\varphi,\xi,\eta)$ is the structure (\ref{fiksieta}) restricted to $M$ and $g$ is the Lorentz metric given by 
$$
 \left[g(\partial_i,\partial_j)\right]=\left[\begin{array}{ccc}
                           -2z &         \ 0 & \  0 \\
                             0 & \ 4x^2 + 2z & \ 2x \\
                             0 &        \ 2x & \  1
                            \end{array}\right].
$$

For the Levi-Civita connection, we have 
\begin{eqnarray*}
  & \nabla_{\partial_1}\partial_1 =- \dfrac{x}{z}\;\partial_2+\Big(1+\dfrac{2x^2}{z}\Big)\partial_3, \quad
       \nabla_{\partial_1}\partial_2 = \nabla_{\partial_2}\partial_1 
       = \dfrac{x}{z}\;\partial_2+\Big(1-\dfrac{2x^2}{z}\Big)\partial_3, &\\
  & \nabla_{\partial_1}\partial_3 = \nabla_{\partial_3}\partial_1 
       = \nabla_{\partial_2}\partial_3 = \nabla_{\partial_3}\partial_2 
       = \dfrac{1}{2z}\;\partial_1+\dfrac{1}{2z}\;\partial_2-\dfrac{x}{z}\;\partial_3, &\\
  & \nabla_{\partial_2}\partial_2 = \dfrac{2x}{z}\;\partial_1+\dfrac{x}{z}\;\partial_2 
       -\Big(1+\dfrac{2x^2}{z}\Big)\partial_3,\quad
       \nabla_{\partial_3}\partial_3 = 0. &
\end{eqnarray*}
Using the above and (\ref{c}), we find $\alpha=\beta=(2z)^{-1}$. 

$$
  \gamma(t)=(0, -2\sqrt{-t},t),\quad t<0 \leqno{(a)}
$$ 
is a Frenet slant curve. For such a curve $\varepsilon_1=-1$, $c=1$, $\alpha(\gamma(t))=\beta(\gamma(t))=\dfrac{1}{2t}$, $\delta(t)=\dfrac{1}{t}$, $\kappa=-\dfrac{\sqrt{5}}{\sqrt{2}t}$, $\tau=-\dfrac{3}{2t}$.


$$
	\gamma(t)=(1/4,t,3/8),  \quad  \leqno{(b)}
$$
is a slant curve with null normal and $\varepsilon_1=1$.

For such a curve, we have $c=1/2$, $\alpha(\gamma(t))=4/3$, $\beta(\gamma(t))=4/3$
and
\begin{eqnarray*}
&& T=(0,1,0),\quad N=(4/3,2/3,-4/3)\\
&& W=(-1/2, 1/4,-1/2) , \quad \kappa=-2/3.
\end{eqnarray*}

$$
	\gamma(t)=(\sqrt{t},-a\sqrt{t},at), \quad t>0, \quad 
	 a=\sqrt[3]{1-b}+\sqrt[3]{1+b},\quad b=\sqrt{\frac{26}{27}}, \leqno{(c)}
$$
is a Legendre curve with null normal and $\varepsilon_1=1$.
For such a curve, we have $\alpha(\gamma(t))=(2at)^{-1}$, $\beta(\gamma(t))=(2at)^{-1}$, $\alpha(\gamma(t))=g(N,\varphi\dot\gamma)=1/(2at)$
and
\begin{eqnarray*}
&& T=((2\sqrt{t})^{-1},-a(2\sqrt{t})^{-1},a),\quad N=(1/4t^{-3/2},-1/(4a)t^{-3/2},0),\\
&& W=(-2a^2\sqrt{t},2a\sqrt{t},-8at) , \quad \kappa=\dfrac{1-2a}{2at}.
\end{eqnarray*}

}
\end{ex}
\begin{ex}
{\rm   
Suppose that $M=\mathbb R^2\times\mathbb R_+$. Define a normal almost paracontact structure $(\varphi,\xi,\eta)$ on $M$ by 
\begin{equation}
  \varphi\partial_1=\partial_2,\quad
  \varphi\partial_2=\partial_1,\quad
  \varphi\partial_3=0,\quad
  \xi=\partial_3,\quad
  \eta=dz,
\end{equation}
and compatible with this structure a Lorentz metric  
$$
 \left[g(\partial_i,\partial_j)\right]=\left[\begin{array}{ccc}
                           -2z & \ 0 & \  0 \\
                             0 & \  2z & \ 0 \\
                             0 &  \ 0 & \  1
                            \end{array}\right].
$$
The quadruple $(\varphi,\xi,\eta,g)$ becomes a normal almost paracontact metric structure on $M$. For the Levi-Civita connection, we find
\begin{eqnarray*}
  && \nabla_{\partial_1}\partial_1 = \partial_3, \quad
     \nabla_{\partial_1}\partial_2 = \nabla_{\partial_2}\partial_1 = 0,\quad
     \nabla_{\partial_1}\partial_3 = \nabla_{\partial_3}\partial_1 = \frac{1}{2z}\;\partial_1,\\
  && \nabla_{\partial_2}\partial_2 = -\partial_3,\quad
     \nabla_{\partial_2}\partial_3 = \nabla_{\partial_3}\partial_2 = \frac{1}{2z}\;\partial_2, \quad
     \nabla_{\partial_3}\partial_3 = 0.
\end{eqnarray*}
Using the above and (\ref{c}), we get $\alpha=(2z)^{-1}$ and $\beta=0$. 

The curve
$$
 \gamma(t)=\Big(\cosh t,\sinh t,\frac{1}{2}\Big),
$$
is a Legendre curve with $\varepsilon_1=1$ and null normal in $M$.

Then for $\gamma$, we get $\alpha(\gamma(t))=1$, $\delta(t)=-1$ and
\begin{eqnarray*}
&&  T=(\sinh t,\cosh t,0),\quad N=\left(\cosh t, \sinh t, -1\right),\\
&&  W= -\dfrac{1}{2}\left(\cosh t, \sinh t, 1\right), \quad \kappa=0.
\end{eqnarray*}

}
\end{ex}


\end {document}